\documentclass [11pt, leqno]{amsart}

\begin{document}

\begin{center}
  {\large\bf \uppercase{On a three-dimensional Riemannian\\[3pt] manifold with an additional structure}}

\vskip1pc
  {\bf Georgi Dzhelepov, Iva Dokuzova, Dimitar Razpopov}
\end{center}

\vskip.5cm

\begin{quote}
  \hskip 0.6cm \small {\bf Abstract.} We consider a $3$-dimensional Riemannian manifold $M$ with a
metric tensor $g$, and affinors $q$ and
$S$. We note that the local coordinates of these three tensors are
circulant matrices. We have that the third degree of $q$ is the identity and $q$ is compatible with $g$. We discuss the sectional curvatures in case when $q$ is parallel with respect to the connection of $g$.
\end{quote}



\thispagestyle{empty}
\bigskip

\bigskip
\centerline{\bf 1. Introduction}
\bigskip

 Many papers in the differential geometry have been dedicated on the problems in the differential manifolds admitting an additional affinor structure~$f$. In the most of them $f$ satisfies some identities of the second degree $f^{2}=id$, or $f^{2}=-id$. We note two papers [7], [8] where $f$ satisfies the equation of the third degree $f^{3}+f=0$.

Let a differential manifold admit an affine connection $\nabla$ and an affinor structure $f$. If $\nabla f$ satisfies some equation there follows an useful   curvature identity. Such identities and assertions were obtained in the almost Hermitian geometry in [2]. Analogous results have been discussed for the almost complex manifolds with Norden metric in [1], [3] and [4], and for the almost contact manifolds with $B$-metric in [5] and [6].

 In the present paper we are interested in a three-dimensional Riemannian manifold $M$ with an affinor structure $q$. The structure satisfies the identity $q^{3}=id$, $q\neq \pm id$ and $q$ is compatible with the Riemannian metric of $M$. Moreover, we suppose the local coordinates of these structures are circulant. We search conditions the structure $q$ to be parallel with respect to the Riemannian connection $\nabla$ of
$g$ (i.e. $\nabla q =0$). We get some curvature identities in this case.

\bigskip\bigskip
\centerline{\bf 2. Preliminaries}
\bigskip

It is known from the linear algebra, that the set of circulant matrices of type
$(n\times n)$ is a commutative group. In the present paper we use four
circulant matrices of type $(3\times 3)$ for geometrical
considerations, as follows:
\begin{equation}\label{f2}
    \big(g_{ij}\big)=\begin{pmatrix}
      A & B & B \\
      B & A & B \\
      B & B & A \\
    \end{pmatrix}, \quad A>B>0,
\end{equation}
where $A=A\big(X^{1},\, X^{2},\, X^{3}\big)$, $B=B\big(X^{1},\, X^{2},\, X^{3}\big)$; and $X^{1},\, X^{2},\, X^{3}\in R$.
\begin{equation}\label{f3}
    \big(g^{ij}\big)=\frac{1}{D}\begin{pmatrix}
      A+B & -B & -B \\
      -B & A+B & -B \\
      -B & -B & A+B \\
    \end{pmatrix}, \quad D=\big(A-B\big)\big(A+2B\big),
\end{equation}
\begin{equation}\label{f4}
    \left(q_{i}^{.j}\right)=\begin{pmatrix}
      0 & 1 & 0 \\
      0 & 0 & 1 \\
      1 & 0 & 0 \\
    \end{pmatrix},
\end{equation}
\begin{equation}\label{f5}
    \left(S_{i}^{.j}\right)=\begin{pmatrix}
      -1 & 1 & 1 \\
      1 & -1 & 1 \\
      1 & 1 & -1 \\
    \end{pmatrix}.
\end{equation}

We choose the form in (\ref{f4}) of the matrix $q$ because of the next assertion:

\bigskip

{\bf Lemma 1.} {\it  Let $\big(m_{ij}\big)$, $i,\,j=1,\,2,\,3$ be a circulant non-degenerate matrix and its third degree is the unit matrix.

Then $\big(m_{ij}\big)$ has one of the following forms: }
\begin{equation}\label{m}
\begin{pmatrix}
      1 & 0 & 0 \\
      0 & 1 & 0 \\
      0 & 0 & 1 \\
    \end{pmatrix}, \hspace{0.5 cm}
     \begin{pmatrix}
      0 & 1 & 0 \\
      0 & 0 & 1 \\
      1 & 0 & 0 \\
    \end{pmatrix}, \hspace{0.5 cm} \begin{pmatrix}
      0 & 0 & 1 \\
      1 & 0 & 0 \\
      0 & 1 & 0 \\
    \end{pmatrix}.
\end{equation}

{\it Proof.} If $\big(m_{ij}\big)$ has the form
\begin{equation*}
    \big(m_{ij}\big)=\begin{pmatrix}
      a & b & c \\
      c & a & b \\
      b & c & a \\
    \end{pmatrix},
\end{equation*}
then from the condition $\big(m_{ij}\big)^{3}=E$ ($E$ is the unit matrix) we get the system
\begin{align*}
    &a^{3}+b^{3}+c^{3}+6\,a\,b\,c=1\\
    &a^{2}b+ac^{2}+b^{2}c=0\\
    &ab^{2}+ca^{2}+c^{2}b=0.
\end{align*}

The all solutions of this system are (\ref{m}).

%%% ----------------------------------------------------------------

\bigskip\bigskip
\centerline{\bf 3. A Parallel Structure}
\bigskip

Let $M$ be a $3$-dimensional Riemannian manifold and $\big\{e_{1},\, e_{2},\, e_{3}\big\}$ be a basis of the tangent space $T_{p}M$ at a point $p\big(X^{1},\, X^{2},\, X^{3}\big)\in M$. Let $g$ be a metric tensor and $q$ be an affinor, which local coordinates are given in (\ref{f2})~and~(\ref{f4}), respectively.
Let $A$ and $B$ from (\ref{f2}) be smooth functions of a point
$p$ in some coordinate neighborhood $F\subset R^{3}$. We will use the
notation $\Phi_{i}=\dfrac{\partial \Phi}{\partial X^{i}}$ for
every smooth function $\Phi$, defined in $F$. We verify that the following identities are true
\begin{equation}\label{2.1}
    q^{3}=E;\quad g(qx, qy)=g(x,y),\quad x,\ y\in \chi M,
\end{equation}
 as well as
\begin{equation}\label{2.2}
    g_{is}\,g^{js}=\delta_{i}^{j}.
\end{equation}

Let $\nabla$ be the Riemannian connection of $g$ and
$\Gamma_{ij}^{s}$ be the Christoffel symbols of $\nabla$. It is
well known the next formula
\begin{equation}\label{2.3}
2\Gamma_{ij}^{s}=g^{as} \left(\partial_{i}g_{aj}+\partial_{j}g_{ai}-\partial_{a}g_{ij}\right).
\end{equation}

Using (\ref{f2}), (\ref{f3}), (\ref{2.2}), (\ref{2.3}), after long computations we get the next equalities:
\begin{align}\label{2.4}\nonumber
\Gamma_{ii}^{i}&=\frac{1}{2D}\left((A+B)A_{i}-B(4B_{i}-A_{j}-A_{k})\right),\\\nonumber
\Gamma_{ii}^{k}&=\frac{1}{2D}\left((A+B)(2B_{i}-A_{k})-B(2B_{i}-A_{j}+A_{i})\right),\\
\Gamma_{ij}^{i}&=\frac{1}{2D}\left((A+B)A_{j}-B(-B_{k}+B_{i}+B_{j}+A_{i})\right),\\\nonumber
\Gamma_{ij}^{k}&=\frac{1}{2D}\left((A+B)(-B_{k}+B_{i}+B_{j})-B(A_{i}+A_{j})\right),
\end{align}
where $i\neq j\neq k$ and $i=1,\, 2,\, 3$, $j=1,\, 2,\, 3$, $k=1,\, 2,\, 3$.

   \bigskip

{\bf Theorem 1.}  {\it  Let $M$ be the Riemannian manifold, supplied with a metric tensor $g$, and affinors  $q$ and $S$, defined by (\ref{f2}), (\ref{f4}) and (\ref{f5}), respectively.  The structure $q$ is parallel with respect to the Riemannian connection $\nabla$ of $g$,
    if and only if}
    \begin{equation}\label{2.5}
    \textit{grad}\, A=\textit{grad}\, B . S.
\end{equation}

{\it Proof.}

\noindent a) Let $q$ be a parallel structure with respect to $\nabla$, i.e.
\begin{equation}\label{dop2.1}
\nabla
q=0.
\end{equation}

In terms of the local coordinates, the last equation implies
$$\nabla_{i}q_{j}^{s}=\partial_{i}q_{j}^{s}+\Gamma_{ia}^{s}q_{j}^{a}-\Gamma_{ij}^{a}q_{a}^{s}=0,$$
 which, by virtue of (\ref{f4}), is equivalent to
\begin{equation}\label{2.6}
    \Gamma_{ia}^{s}q_{j}^{.a}=\Gamma_{ij}^{a}q_{a}^{.s}.
\end{equation}

Using (\ref{f4}), (\ref{2.4}) and (\ref{2.6}), we get $18$
equations which all imply (\ref{2.5}).

\bigskip

\noindent b) Vice versa, let (\ref{2.5}) be valid. Then from (\ref{2.4}) we get
\begin{align*}
    \Gamma_{11}^{1}=\Gamma_{12}^{2}=\Gamma_{13}^{3}=\Gamma_{22}^{3}=\Gamma_{23}^{1}=\Gamma_{33}^{2}&=\frac{1}{2D}(AA_{1}+B(-3B_{1}+B_{2}+B_{3})),\\ \Gamma_{11}^{3}=\Gamma_{12}^{1}=\Gamma_{13}^{2}=\Gamma_{22}^{2}=\Gamma_{23}^{3}=\Gamma_{33}^{1}&=\frac{1}{2D}(AA_{2}+B(B_{1}-3B_{2}+B_{3})),\\
    \Gamma_{11}^{2}=\Gamma_{12}^{3}=\Gamma_{13}^{1}=\Gamma_{22}^{1}=\Gamma_{23}^{2}=\Gamma_{33}^{3}&=\frac{1}{2D}(AA_{3}+B(B_{1}+B_{2}-3B_{3})).
    \end{align*}

Now, we can verify that (\ref{2.6}) is valid. That means
$\nabla_{i}q_{j}^{s}=0$, i.e. $\nabla q=0$.

\hfill$\square$

\bigskip

{\bf Remark.} {\it In fact (\ref{2.5}) is a system of three partial
differential equations for the functions $A$ and $B$. Let
$p(X^{1}, X^{2}, X^{3})$ be a point in $M$.  We assume
$B=B(p)$ as a known function and then we can say that (\ref{2.5}) has a solution.
Particularly, we give a simple (but non-trivial example) for both
functions, satisfying (\ref{2.5}), as follows $A=(X^{1})^{2}+(X^{2})^{2}+(X^{3})^{2}$; $B=X^{1}X^{2}+X^{1}X^{3}+X^{2}X^{3}$, where $A > B > 0$. }

\bigskip\bigskip
\centerline{\bf 4. Sectional Curvatures}
\bigskip

Let $M$ be the Riemannian manifold with a metric tensor $g$ and a structure~$q$, defined by (\ref{f2}) and (\ref{f4}),
respectively.
Let $R$ be the curvature tensor field of $\nabla$, i.e $R(x, y)z=\nabla_{x}\nabla_{y}z-\nabla_{y}\nabla_{x}z-\nabla_{[x,y]}z$. We consider the associated tensor field $R$ of type $(0, 4)$, defined by the condition
\begin{equation*}
    R(x, y, z, u)=g\big(R(x, y)z,u\big), \qquad x,\, y,\, z,\, u\in \chi M.
    \end{equation*}

{\bf Theorem 2.} {\it If $M$ is the Riemannian manifold with a metric tensor $g$ and a
parallel structure $q$, defined by (\ref{f2}) and (\ref{f4}),
respectively, then the curvature tensor $R$ of $g$ satisfies the
identity:}
\begin{equation}\label{3.1}
    R(x, y, q^{2}z, u)=R(x, y, z, qu),\qquad x, y, z, u\in \chi M.
\end{equation}

{\it Proof.} In terms of the local coordinates (\ref{dop2.1}) implies
\begin{equation}\label{3.2}
    R^{l}_{sji}\,q_{k}^{.s}=R^{s}_{kji}\,q_{s}^{.l}.
\end{equation}

Using (\ref{f4}), we verify $q^{i}_{.j}=q_{a}^{.i}\,q_{j}^{.a}$
and then from (\ref{f2}), (\ref{f3}) and (\ref{3.2}) we obtain
(\ref{3.1}).

\hfill$\square$

\smallskip

Let $p$ be a point in $M$ and $x$, $y$ be two linearly
independent vectors on~$T_{p}M$. It is known that the quantity
\begin{equation}\label{3.3}
    \mu(\textsl{L};p)=\frac{R(x, y, x, y)}{g(x, x)g(y, y)-g^{2}(x, y)}
\end{equation}
is the sectional curvature of $2$-plane $\textsl{L}=\{x, y\}$.

Let $p$ be a point in $M$ and $x=(x^{1},\, x^{2},\, x^{3})$ be a vector in $T_{p}M$. The vectors $x,\ qx,\ q^{2}x$ are linearly independent, when
   \begin{equation}\label{lema2}
   3x^{1}x^{2}x^{3}\neq (x^{1})^{3}+(x^{2})^{3}+(x^{3})^{3}.
   \end{equation}

    Then we define
$2$-planes $\textsl{L}_{1}=\{x, qx\}$, $\textsl{L}_{2}=\{qx, q^{2}x\}$ and $\textsl{L}_{3}=\{q^{2}x, x\}$ and we prove the following

\bigskip

{\bf Theorem 3.} {\it
Let $M$ be the Riemannian manifold with a metric tensor $g$ and a
parallel structure $q$, defined by (\ref{f2}) and (\ref{f4}),
respectively. Let $p$ be a point in $M$ and $x$ be an
arbitrary vector in $T_{p}M$ satisfying (\ref{lema2}). Then the sectional curvatures of
$2$-planes $\textsl{L}_{1}=\{x, qx\}$, $\textsl{L}_{2}=\{qx, q^{2}x\}$,
$\textsl{L}_{3}=\{q^{2}x, x\}$ are equal. }

{\it Proof.} From (\ref{3.1}) we obtain
\begin{equation}\label{3.6}
    R(x,\, y,\, z,\, u)=R(x,\, y,\, qz,\, qu)=R(x,\, y,\, q^{2}z,\, q^{2}u).
    \end{equation}

    In (\ref{3.6}) we set the following substitutions:
    a) $z=x$, $y=u=qx$; b)~$x\sim qx$, $z=qx$, $y=u=q^{2}x$; c) $x\sim q^{2}x$, $z=q^{2}x$, $y=u=x$.
    Comparing the obtained results, we get
    \begin{equation}\label{dopl}
    \begin{split}
    R(x,\, qx,\,q^{2}x,\, x) &=R(x,\, qx,\, qx,\, q^{2}x)\\
    &=R(q^{2}x,\, x,\, qx,\, q^{2}x) \\
    &=R(x,\, qx,\, x,\, qx)
    \end{split}
    \end{equation}
    and
    \begin{equation}\label{3.7}
    R(x,\, qx,\, x,\, qx)=R(qx,\, q^{2}x,\, qx,\, q^{2}x)=R(q^{2}x,\, x,\, q^{2}x,\, x).
    \end{equation}

Equalities (\ref{2.1}), (\ref{3.3}), (\ref{lema2}) and (\ref{3.7}) imply
    \begin{equation*}
    \mu(\textsl{L}_{1};p)= \mu(\textsl{L}_{2};p)= \mu(\textsl{L}_{3};p)=\frac{R(x,\, qx,\, x, \, qx)}{g^{2}(x, x)-g^{2}(x, qx)}\,.
\end{equation*}

 By virtue of the linear independence of
$x$ and $qx$, we have
$$g^{2}(x,x)-g^{2}(x, qx)=g^{2}(x,x)(1-\cos\varphi)\neq 0, $$ where
$\varphi$ is the angle between $x$ and $qx$.

\hfill$\square$

\bigskip\bigskip
\centerline{\bf 5. An Orthonormal $q$-Base of Vectors in $T_{p}M$}
\bigskip

Let $M$ be the Riemannian manifold with a metric tensor $g$ and a structure~$q$, defined by (\ref{f2}) and (\ref{f4}),
respectively. We note that the only real eigenvalue and the only eigenvector of the structure $q$ are $\lambda=1$ and $x\big(x^{1},\, x^{1},\, x^{1}\big)$, respectively.

Now, let \begin{equation}\label{n2.3}
x=\big(x^{1}, \,x^{2},\, x^{3}\big)
   \end{equation}
    be a non-eigenvector vector of the structure $q$.
We have
\begin{equation}\label{n2.4}
    g(x,x)=\|x\|\|x\| \cos 0=\|x\|^{2}, \quad g(x,qx)=\|x\|\|qx\| \cos\varphi= \|x\|^{2} \cos\varphi,
\end{equation}
where $\|x\|$ and $\|qx\|$ are the norms of $x$ and $qx$; and $\varphi$
is the angle between $x$ and $qx$.

From (\ref{f2}), (\ref{n2.3}) and (\ref{n2.4}) we
calculate
\begin{equation}\label{n2.7}
    g(x,x)=A\left((x^{1})^{2}+(x^{2})^{2}+(x^{3})^{2}\right)+2B\left(x^{1}x^{2}+x^{1}x^{3}+x^{2}x^{3}\right),
\end{equation}
\begin{equation}\label{n2.7*}
    g(x,qx)=B\left((x^{1})^{2}+(x^{2})^{2}+(x^{3})^{2}\right)+(A+B)\left(x^{1}x^{2}+x^{1}x^{3}+x^{2}x^{3}\right).
\end{equation}

The above equations imply $\|x\|=\|qx\|>0$.

\bigskip

{\bf Theorem 4.}{\it Let $M$ be the Riemannian manifold with a metric tensor $g$ and an
affinor structure $q$, defined by (\ref{f2}) and (\ref{f4}),
respectively. Let $x(x^{1}, x^{2}, x^{3})$ be a non-eigenvector on
$T_{p}M$. If $\varphi$ is the angle between $x$ and $qx$, then
we have $\varphi\in \left(0,\dfrac{2\pi}{3}\right)$. }

{\it Proof.} We apply equations (\ref{n2.7}) and (\ref{n2.7*}) in
$\cos\varphi=\dfrac{g(x, qx)}{g(x,x)}\,,$ and we get
\begin{equation}\label{123}  \cos\varphi=\frac{\big((x^{1})^{2}+(x^{2})^{2}+(x^{3})^{2}\big)+(A+B)\big(x^{1}x^{2}+x^{1}x^{3}+x^{2}x^{3}\big)}{A\big((x^{1})^{2}+(x^{2})^{2}+(x^{3})^{2}\big)+2B\big(x^{1}x^{2}+x^{1}x^{3}+x^{2}x^{3}\big)}.
\end{equation}

 Also we have $x(x^{1},\, x^{2},\, x^{3})\neq (x^{1},\, x^{1},\, x^{1})$ because $x$ is a non-eigenvector of~$q$.

We suppose that $\varphi\geq\dfrac{2\pi}{3}$, i.e. $\cos\varphi
\leq-\dfrac{1}{2}$\,. The last condition and (\ref{123}) imply
\begin{equation*}
\frac{B\big((x^{1})^{2}+(x^{2})^{2}+(x^{3})^{2}\big)+(A+B)\big(x^{1}x^{2}+x^{1}x^{3}+x^{2}x^{3}\big)}{A\big((x^{1})^{2}+(x^{2})^{2}+(x^{3})^{2}\big)+2B\big(x^{1}x^{2}+x^{1}x^{3}+x^{2}x^{3}\big)}\leq -\frac{1}{2}\,
\end{equation*}
that gives the inequality
$$(2B+A)\Big((x^{1})^{2}+(x^{2})^{2}+(x^{3})^{2}+2\big(x^{1}x^{2}+x^{1}x^{3}+x^{2}x^{3}\big)\Big)\leq 0.$$

From the condition $A+2B>0$ we get that
$$(x^{1})^{2}+(x^{2})^{2}+(x^{3})^{2}+2\big(x^{1}x^{2}+x^{1}x^{3}+x^{2}x^{3}\big)\leq 0$$
and $(x^{1}+ x^{2}+ x^{3})^{2}\leq 0$. The last inequality has no
solution in the real set. Then we have $\cos\varphi
>-\dfrac{1}{2}$.

\hfill$\square$

\smallskip

Immediately, from Theorem~4, we establish that an orthonormal $q$-base $(x,\, qx,\, q^{2}x)$ in $T_{p}M$ exists.
Particularly, we verify that the vector
\begin{equation}\label{exampl}
x=\left(\frac{\sqrt{A-B}+\sqrt{A+3B}}{2\sqrt{A^{2}+AB-2B^{2}}}\,,\quad
\frac{\sqrt{A-B}-\sqrt{A+3B}}{2\sqrt{A^{2}+AB-2B^{2}}}\,,\quad 0\right)
\end{equation}
satisfies the conditions
\begin{equation}\label{n3.1}
    g(x, x)=1, \qquad  g(x, qx)=0.
\end{equation}

The base $(x,\, qx,\, q^{2}x)$, where $x$ satisfies (\ref{exampl}), is an example of an orthonormal $q$-base in $T_{p}M$.

\bigskip

{\bf Theorem 5} {\it
Let $M$ be the Riemannian manifold with a metric tensor $g$ and a
parallel structure $q$, defined by (\ref{f2}) and (\ref{f4}),
respectively. Let $(x,\, qx,\, q^{2}x)$ be an orthonormal $q$-base in
$T_{p}M$, $p\in M,$ and $u=\alpha\,.\,x+\beta\, .\,qx +\gamma
\,.\,q^{2}x$,  $v=\delta .\,x+\zeta.\,qx +\eta .\,q^{2}x$ be arbitrary
vectors in $T_{p}M$. For the sectional curvature $\mu(u, v)$
of $2$-plane $\{u, v\}$ we have }
\begin{equation}\label{cor}
    \mu(u,v)=\frac{\big(\alpha\zeta-\beta\delta+\delta\gamma-\alpha\eta +\beta\eta-\gamma\zeta\big)^{2}}{\big(\alpha^{2}+\beta^{2}+\gamma^{2}\big)\big( \delta^{2}+\zeta^{2}+\eta^{2}\big)-\big(\alpha\delta+\beta\zeta+\gamma\eta\big)^{2}}\,\mu(x, qx).
\end{equation}

{\it Proof.}  We calculate
\begin{equation}\label{n4.1}
    g(u, u)= \alpha^{2}+\beta^{2}+\gamma^{2},\quad g(v, v)= \delta^{2}+\zeta^{2}+\eta^{2},
\end{equation}
  $$g(u, v)= \alpha\delta+\beta\zeta+\gamma\eta .$$

For the sectional curvature of $2$-plane $\{u, v\}$ we have
\begin{equation}\label{n4.2}
    \mu(u,v)=\frac{R(u,\, v,\, u,\, v)}{g(u, u)g(v, v)-g^{2}(u, v)}\,.
\end{equation}

 Using the linear properties of the metric $g$ and the curvature tensor field~$R$ after long calculations we get
\begin{equation}\label{n4.3}\begin{split}
    R(u,\,v,\,u,\,v)= &(\alpha\zeta-\beta\delta)^{2}R(x, \,qx,\, x,\, qx)\\
    &\hspace{0.4 cm}+(\delta\gamma-\alpha\eta)^{2}R(x,\, q^{2}x,\, x, \,q^{2}x)\\
    & \hspace{0.5 cm} +(\beta\eta-\gamma\zeta)^{2}R(qx,\, q^{2}x,\, qx,\, q^{2}x)\\
    & \hspace{0.6 cm} +2(\alpha\zeta-\beta\delta)(\delta\gamma-\alpha\eta)R(x,\, qx, \,q^{2}x,\, x)\\
    & \hspace{0.7 cm} +2(\delta\gamma-\alpha\eta)(\beta\eta-\gamma\zeta)R(q^{2}x,\, x,\, qx,\, q^{2}x)\\
    & \hspace{0.8 cm}+2(\alpha\zeta-\beta\delta)(\beta\eta-\gamma\zeta)R(x,\, qx,\, qx,\, q^{2}x).\\
    \end{split}
\end{equation}

    From (\ref{dopl}), (\ref{3.7}) and (\ref{n4.3}) we obtain
    \begin{equation}\label{n4.6}
    R(u,\,v,\,u,\,v)= \big((\alpha\zeta-\beta\delta) +(\delta\gamma-\alpha\eta) +(\beta\eta-\gamma\zeta)\big)^{2}R(x,\, qx,\, x,\, qx).
    \end{equation}

    From (\ref{n4.1}), (\ref{n4.2}) and (\ref{n4.6}) we get
    \begin{equation*}
    \mu(u,v)=\frac{\big(\alpha\zeta-\beta\delta +\delta\gamma-\alpha\eta +\beta\eta-\gamma\zeta\big)^{2}}{(\alpha^{2}+\beta^{2}+\gamma^{2})( \delta^{2}+\zeta^{2}+\eta^{2})-(\alpha\delta+\beta\zeta+\gamma\eta)^{2}}R(x, \,qx,\, x,\, qx).
\end{equation*}

The last equation and (\ref{n3.1}) imply (\ref{cor}).

\hfill$\square$

\bigskip

{\bf Corollary 1.} {\it  Let $u$ be an arbitrary non-eigenvector in $T_{p}M$, $p\in
M$, and~$\theta$ be the angle between $u$ and $qu$.

Then we have }
\begin{equation}\label{mju}
    \mu(u,qu)=\mu(x, qx)\tan^{2}\frac{\theta}{2}\,,\quad \theta\in \left(0,\frac{2\pi}{3}\right).
\end{equation}

{\it Proof.} In (\ref{cor}) we substitute $v=qu$, $\delta=\gamma$, $\zeta=\alpha$, $\eta=\beta$ and we obtain
\begin{equation*}
    \mu(u,qu)=\frac{(\alpha^{2}+\beta^{2}+\gamma^{2}-\beta\gamma
    -\alpha\beta
    -\alpha\gamma)^{2}}{(\alpha^{2}+\beta^{2}+\gamma^{2})^{2}-(\alpha\gamma+\alpha\beta+\gamma\beta)^{2}}\,\mu(x, qx).
\end{equation*}

Then from (\ref{n4.1}) we get
$$\mu(u,qu)=\dfrac{(g(u, u)-g(u, qu))^{2}}{g^{2}(u,u)-g^{2}(u, qu)}\,\mu(x, qx),$$
i.e.
$$\mu(u,qu)=\dfrac{(1-\cos\theta)^{2}}{1-\cos^{2}\theta}\,\mu(x, qx),$$
which implies (\ref{mju}).

\hfill$\square$

\bigskip

{\bf Corollary 2.} {\it Let $u, v$ be an arbitrary non-eigenvectors on $T_{p}M$, $p\in M$, and~$\theta$ be the angle between $u$ and $qu$, and $\psi$ be
the angle between $v$ and $qv$.

Then we have }
\begin{equation*}
    \mu(u,qu)\tan^{2}\frac{\psi}{2}=\mu(v,
    qv)\tan^{2}\frac{\theta}{2},\quad \psi,\, \theta \in \left(0,\frac{2\pi}{3}\right).
\end{equation*}

The proof follows immediately from (\ref{mju}).

\bigskip\bigskip
\centerline{\bf Acknowledgments}
\bigskip

Research was partially supported by  project RS11 - FMI - 004 of the Scientific Research Fund, Paisii Hilendarski University of Plovdiv, Bulgaria.

\bigskip
%\newpage
\begin{center}
{\large\bf References}
\end{center}

\smallskip

\leftskip 1.2pc
\parindent-1.2pc

[1] {\sc{Borisov A., Ganchev G.,}}  Curvature properties of Kaehlerian manifolds with B-metric, {\it Math. Educ. Math., Proc of 14th Spring Conf. of UBM},
Sunny Beach, (1985), 220--226.

[2] {\sc  Grey A.,}  Curvature identities for Hermitian and Almost Hermitian Manifolds, {\it Tohoku Math. Journal}, Vol. 28, No. 4, (1976), 601--612.

[3] {\sc  Gribachev K., Mekerov  D., Djelepov G.,} On the Geometry of Almost $B$-manifolds, {\it Compt. Rend. Acad. Bulg. Sci.}, Vol. 38, No. 5, (1985), 563--566.

[4] {\sc Gribachev  K., Djelepov G.,} On the Geometry of the normal generalized $B$-manifolds, {\it PU Sci. Works-math.}, Vol. 23, No. 1,  (1985), 157--168.

[5] {\sc Manev M.,  Nakova G.,}  Curvature properties on some three-dimensional
almost contact B-metric manifolds, {\it Plovdiv Univ. Sci. Works - math.}, Vol. 34, no. 3, (2004), 51–-60.

[6] {\sc Nakova G., Manev M.,}  Curvature properties on some three-dimensional almost contact manifolds with B-metric,{\it Proc.
5th Int. Conf. Geometry, Integrability \& Quantization V Eds.I. M. Mladenov, A. C. Hirshfeld}, SOFTEX, Sofia, (2004), 169–-177.

[7] {\sc Yano K.,  Ishihara S.,}  Structure defined by $f$  $f^{3}+f=0$, {\it Proc.US-Japan Seminar in Differential Geometry.}, Kyoto, (1965), 153--166.

[8] {\sc Yano K.,}  On a structure defined by a tensor field of type $(1, 1)$ satisfying $f^{3}+f=0$, Tensor, (1963), 99--109.

\leftskip 0pc
\parindent .7cm
\vskip1pc\vskip3pt\medskip

\bigskip

\noindent
Georgi Dzhelepov\\
Department of Mathematics and Physics \\
Agricultural University of Plovdiv\\
12 Mendeleev Blvd., 4000 Plovdiv, Bulgaria\\
e-mail:  \verb"dzhelepov@au-plovdiv.bg"
\\
\\
Iva Dokuzova\\
Faculty of Mathematics and Informatics\\
University of Plovdiv\\
236 Bulgaria Blvd., 4003 Plovdiv, Bulgaria\\
e-mail: \verb"dokuzova@uni-plovdiv.bg"
\\
\\
Dimitar Razpopov \\
Department of Mathematics and Physics \\
Agricultural University of Plovdiv\\
12 Mendeleev Blvd., 4000 Plovdiv, Bulgaria\\
e-mail: \verb"razpopov@au-plovdiv.bg"
\\
\\

\bigskip
\bigskip

\end{document}